\documentclass[12pt]{article}

\usepackage{amsmath}
\usepackage{amssymb}
\usepackage{amsfonts}
\usepackage{amscd}
\usepackage{a4}
\usepackage{graphicx}

\allowdisplaybreaks
\numberwithin{equation}{section}
\newtheorem{Theorem}{Theorem}[section]
\newtheorem{Corollary}[Theorem]{Corollary}
\newtheorem{Lemma}[Theorem]{Lemma}
\newtheorem{Proposition}[Theorem]{Proposition}

\newenvironment{Remark}{\refstepcounter{Theorem}
\noindent {\it Remark \arabic{section}.\arabic{Theorem}.}}\noindent

\noindent

\newcommand*{\NN}{{\mathbb N}}
\newcommand*{\RR}{{\mathbb R}}
\newcommand*{\ZZ}{{\mathbb Z}}

\newcommand*{\CC}{{\mathbb C}}

\newcommand*{\B}{{\mathcal B}}
\newcommand*{\D}{{\mathcal D}}
\newcommand*{\F}{{\mathcal F}}
\newcommand*{\G}{{\mathcal G}}
\newcommand*{\HH}{{\mathcal H}}

\newcommand*{\LL}{{\mathcal L}}

\newcommand*{\Var}{\text{\rm Var}}

\newcommand*{\II}{I\hspace{-2pt}I}
\newcommand*{\III}{I\hspace{-2pt}I\hspace{-1.9pt}I}

\def\defeq{\mathrel{\mathop:}=}
\def\eqdef{\mathrel{=\hspace{-3pt}\mathop:}}

\begin{document}

\renewcommand{\thefootnote}{}

\begin{center}\Large\bf 

The nonexplosive solution\\
of explosive autoregressions\\
\end{center}

\begin{center}
\it

Erich H\"ausler\hspace{1.5pt}$^1$\footnote{\hspace{-17.5pt}$^1$\hspace{0.4pt}Mathematical Institute, 
University of Giessen, Giessen, Germany, erich.haeusler@t-online.de} and 
Harald Luschgy\hspace{1.5pt}$^2$\footnote{\hspace{-17.5pt}$^2$\hspace{0.5pt}FB IV, Mathematics, University of Trier, Trier, Germany, luschgy@uni-trier.de}

\end{center}\bigskip\bigskip

\centerline{\it Abstract}\bigskip

We establish several features of the stationary solution of purely explosive autoregressions of order $d$ 
based on nonstandard initial values.\bigskip

\section{\large Introduction}
\label{Section:Introduction}

We consider an autoregressive process $\left(Y_n\right)_{n\geq-d+1}$ of order $d$ generated recursively by
\begin{equation}
\label{Eq:Yautoregression}
Y_n=\theta_1Y_{n-1}+\cdots+\theta_dY_{n-d}+Z_n\,,\quad n\geq1\,,
\end{equation}
where $\theta=\left(\theta_1,\ldots,\theta_d\right)^T\in\RR^d,\,\left(Z_n\right)_{n\geq1}$ is an i.i.d.
sequence of real random variables with $Z_1\in\LL^2\left(P\right),\,E\left(Z_1\right)=0,\,\sigma^2=\Var\left(Z_1\right)>0$
and the initial value $\left(Y_0,\ldots,Y_{-d+1}\right)^T$ is an arbitrary $\RR^d$-valued random vector.
Then the process $\left(U_n\right)_{n\geq0}$ given by $U_n=\left(Y_n,Y_{n-1},\ldots,Y_{n-d+1}\right)^T$
can be expressed as a $d$-dimensional autoregression of order one, namely
\begin{equation}
\label{Eq:Unautoregression}
U_n=BU_{n-1}+W_n,\quad n\geq1\,,
\end{equation}
where

\begin{displaymath}
B=B\left(\theta\right)=
\begin{pmatrix}
\theta_1 & \ldots  & \theta_{d-1} & \theta_d \\
         &         &              & 0        \\
         & I_{d-1} &              & \vdots   \\
         &         &              & 0        
\end{pmatrix}
\qquad\text{and}\qquad W_n=\left(Z_n,0,\ldots,0\right)^T=Z_ne_1\,,
\end{displaymath}
with $e_1,\ldots,e_d$ denoting the unit vectors in $\RR^d$. By induction, we have
\begin{equation}
\label{Eq:Un}
U_n=B^nU_0+\sum_{j=1}^nB^{n-j}W_j=B^nU_0+\sum_{j=1}^nB^{j-1}W_{n+1-j}\,,\quad n\geq0\,.
\end{equation}
All random variables are defined on a probability space $\left(\Omega,\F,P\right)$.
\smallskip

Now we concentrate on the purely explosive case. Let $\rho\left(B\right)$ denote the spectral
radius of $B$, i.e. the maximum absolute value of all eigenvalues of $B$ and let 
$\underline{\rho}\left(B\right)$ denote the lower spectral radius of $B$, that is the minimum
absolute value of the eigenvalues of $B$, where the eigenvalues are the solutions in $\CC$ of
\begin{displaymath}
\det\left(B-\lambda I_d\right)=\left(-1\right)^d\left(\lambda^d-\theta_1\lambda^{d-1}-\cdots-
\theta_{d-1}\lambda-\theta_d\right)=0\,.
\end{displaymath}
The purely explosive case corresponds to $\underline{\rho}\left(B\right)>1$. Then $B$ is invertible
or equivalently, $\theta_d\neq0$ since $\det\left(B\right)=\left(-1\right)^{d+1}\theta_d$, and
moreover, $\rho\left(B^{-1}\right)<1$ since the eigenvalues of $B^{-1}$ are $1/\lambda_j, j=1,\ldots,d$,
if $\lambda_j, j=1,\ldots,d$, are the eigenvalues of $B$. Let
$\Theta_{pe}=\left\{\theta\in\RR^d:\underline{\rho}\left(B\left(\theta\right)\right)>1\right\}$ and
$\Theta_{s}=\left\{\theta\in\RR^d:\rho\left(B\left(\theta\right)\right)<1\right\}$. 
\medskip

\begin{Remark}
\label{Remark:d=2}
Let $d=2$. It is well known that the stable case $\rho\left(B\right)<1$ is represented by the open
triangle
\begin{displaymath}
\Theta_s=\left\{\theta\in\RR^2:\left\lvert\theta_1\right\rvert<2,-1<\theta_2<1-\left\lvert\theta_1\right\rvert\right\}\,.
\end{displaymath}
The critical case $\rho\left(B\right)=1$ corresponds to the boundary of $\Theta_s$ so that the parameter
set $\Theta_{pe}$ which corresponds to the purely explosive case $\underline{\rho}\left(B\right)>1$
is (a large) part of the complement of the closed triangle
\begin{displaymath}
\left\{\theta\in\RR^2:\left\lvert\theta_1\right\rvert\leq2,-1\leq\theta_2\leq1-\left\lvert\theta_1\right\rvert\right\}\,.
\end{displaymath}
More precisely, the eigenvalues being
\begin{displaymath}
\lambda_{1,2}=\frac{1}{2}\left(\theta_1\pm\sqrt{\theta_1^2+4\theta_2}\right)\in\CC\,,
\end{displaymath}
one checks that
\begin{displaymath}
\Theta_{pe}=\left\{\theta\in\RR^2:\theta_2>1+\left\lvert\theta_1\right\rvert\right\}\cup
\left\{\theta\in\RR^2:\theta_2<-1,\theta_2<1-\left\lvert\theta_1\right\rvert\right\}\,.
\end{displaymath}
See Figure 1 for the geometric shape of the sets $\Theta_s$ and $\Theta_{pe}$. 
\hfill$\Box$
\end{Remark}
\medskip

\begin{figure}[t]
\textbf{Figure 1.} \textit{The sets} $\Theta_s$ \textit{and} $\Theta_{pe}$ \textit{for} $d=2$
\vspace{-30pt}
\begin{center}
\includegraphics[scale=0.6]{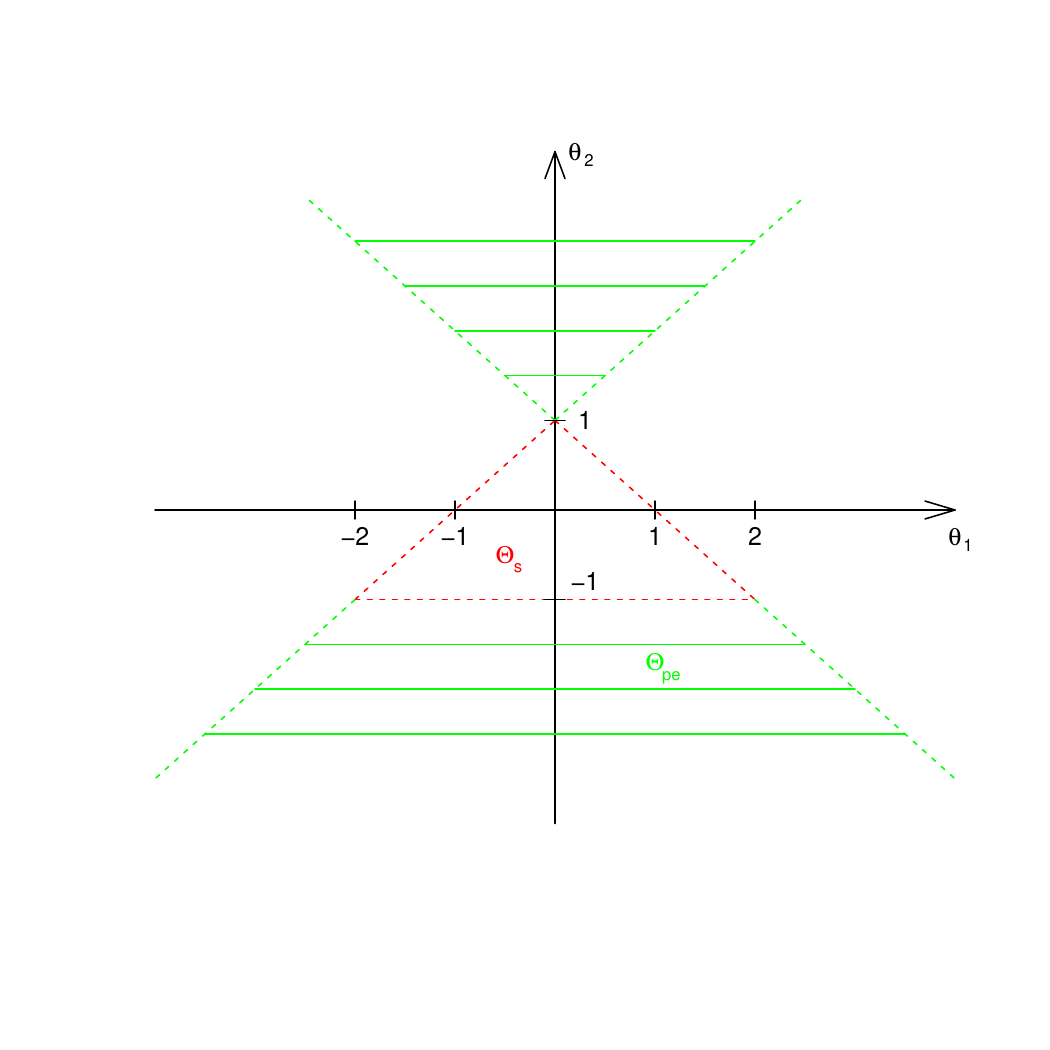}
\end{center}
\vspace{-60pt}
\end{figure}

There exists a submultiplicative norm 
$\left\lVert\cdot\right\rVert_{B^{-1}}$ on $\RR^{d\times d}$ with 
$\left\lVert B^{-1}\right\rVert_{B^{-1}}<1$ (see \cite{Schott}, Theorem~4.24).
Let $\langle\cdot,\cdot\rangle$ denote the standard Euclidean scalar product in $\RR^d$,
let $\left\lVert\cdot\right\rVert$ denote the associated Euclidean $\ell^2$-norm and
let $\left\lVert\cdot\right\rVert_F$ denote the Frobenius norm on $\RR^{d\times d}$
given by $\left\lVert D\right\rVert_F=\left(\sum_{i=1}^d\sum_{j=1}^dD_{ij}^2\right)^{1/2}$
which is submultiplicative, compatible with the $\ell^2$-norm, i.e.
$\left\lVert Dx\right\rVert\leq\left\lVert D\right\rVert_F\left\lVert x\right\rVert$
for every $D\in\RR^{d\times d}$ and $x\in\RR^d$, and satisfies 
$\left\lVert D\right\rVert_F=\left\lVert D^T\right\rVert_F$. Since all norms on
$\RR^{d\times d}$ are equivalent, there is a real constant $C\left(B^{-1}\right)$
such that $\left\lVert D\right\rVert_F\leq C\left(B^{-1}\right)
\left\lVert D\right\rVert_{B^{-1}}$ for all $D\in\RR^{d\times d}$.
\smallskip

By \eqref{Eq:Un}, we have
\begin{equation}
\label{Eq:B-nUnconvergence}
B^{-n}U_n=U_0+\sum_{j=1}^nB^{-j}W_j\rightarrow U_0+\sum_{j=1}^\infty B^{-j}W_j\quad
\text{almost surely as }n\to\infty
\end{equation}
(with $B^{-j}\defeq\left(B^{-1}\right)^j=\left(B^j\right)^{-1}$), where the series converges absolutely
almost surely since
\begin{align*}
\sum_{j=1}^\infty\left\lVert B^{-j}W_j\right\rVert&\leq\sum_{j=1}^\infty\left\lVert B^{-j}\right\rVert_F
\left\lVert W_j\right\rVert=\sum_{j=1}^\infty\left\lVert B^{-j}\right\rVert_F\left\lvert Z_j\right\rvert\\
&\leq C\left(B^{-1}\right)\sum_{j=1}^\infty\left\lVert B^{-j}\right\rVert_{B^{-1}}\left\lvert Z_j\right\rvert
\leq C\left(B^{-1}\right)\sum_{j=1}^\infty\left\lVert B^{-1}\right\rVert_{B^{-1}}^j\left\lvert Z_j\right\rvert
<\infty
\end{align*}
almost surely. This shows the explosive character of $\left(U_n\right)_{n\geq0}$ if 
$U_0+\sum_{j=1}^\infty B^{-j}W_j$ $\neq0$ with positive probability. However, this effect disappears
under the choice
\begin{equation}
\label{Eq:DefinitionU0}
U_0=\left(Y_0,Y_{-1},\ldots,Y_{-d+1}\right)^T\defeq-\sum_{j=1}^\infty B^{-j}W_j\,.
\end{equation}
Then for $n\geq0$ we obtain
\begin{align}
U_n&=B^nU_0+\sum_{j=1}^nB^{n-j}W_j=-\sum_{j=1}^\infty B^{n-j}W_j+\sum_{j=1}^n B^{n-j}W_j\nonumber\\
&=-\sum_{j=n+1}^\infty B^{n-j}W_j=-\sum_{k=1}^\infty B^{-k}W_{n+k}=-\sum_{k=1}^\infty B^{-k}e_1Z_{n+k}\label{Eq:RepresentationUn}
\end{align}
so that $\left(U_n\right)_{n\geq0}$ is stationary and ergodic, since $\left(Z_n\right)_{n\geq1}$ is
stationary and ergodic, and the same holds true for $\left(Y_{n-d+1}\right)_{n\geq0}$ (see the
subsequent Lemma~\ref{Lemma:Stationarity}). The present paper focusses on these stationary and ergodic
sequences $\left(U_n\right)_{n\geq0}$ and $\left(Y_{n-d+1}\right)_{n\geq0}$. By \eqref{Eq:RepresentationUn},
$\left(U_n\right)_{n\geq0}$ is a purely noncausal moving average process w.r.t. the sequence 
$\left(Z_n\right)_{n\geq1}$. Since
$Y_{n-d+1}=\pi_d\left(U_n\right)=-\sum_{k=1}^\infty\left(B^{-k}\right)_{d1}Z_{n+k}$, the same holds true
for $\left(Y_{n-d+1}\right)_{n\geq0}$, where $\pi_j:\RR^d\rightarrow\RR,\pi_j(x)=x_j$, denotes the projection
on the $j$-th coordinate for $1\leq j\leq d$. Such noncausal
autoregressive processes have applications if $n$ is not a linear time parameter, see e.g.~\cite{LiiRosenblatt},
page 6. Even earlier work on noncausal autoregressions is contained in~\cite{Breidt}. 
More recently, they also found applications in the analysis of financial time series,
where now $n$ is quite naturally a time parameter; see e.g.~\cite{DavisSong}, \cite{GourierouxZakoian},
\cite{LanneSaikkonen}. 
\medskip

In Section~\ref{Section:Features}
we will present probabilistic properties of $U_n$ and related processes including a 
mixing central limit theorem for $U_n$ and $Y_n$. In Section~\ref{Section:CLT} we will establish
a central limit theorem for processes of the type $h\left(U_n\right)Z_n$, where $h:\RR^d\rightarrow\RR^m$
denotes any Borel measurable map with $\left\lVert h\left(U_0\right)\right\rVert\in\LL^2\left(P\right)$.
The main tool is a (stable) central limit theorem for reverse multivariate martingale difference sequences 
(which should not be confused with sequences of differences of reverse martingales). As an application
we derive some asymptotic features of the classical least squares estimator of $\theta\in\Theta_{pe}$ 
when the parameter $\theta$ is unknown.
\bigskip

\section{\Large Features of the stationary solution}
\label{Section:Features}

We first confirm the (strict) stationarity and ergodicity of $\left(U_n\right)_{n\geq0}$ 
as defined by~\eqref{Eq:RepresentationUn} and several related processes. Notice that $\left\lVert U_0\right\rVert\in
\LL^2\left(P\right)$ since $\left\lVert U_0\right\rVert\leq\sum_{j=1}^\infty\left\lVert B^{-j}\right\rVert_F
\left\lvert Z_j\right\vert$ so that
\begin{align*}
E\left(\left\lVert U_0\right\rVert^2\right)&\leq\sum_{j=1}^\infty\sum_{k=1}^\infty\left\lVert B^{-j}\right\rVert_F
\left\lVert B^{-k}\right\rVert_FE\left(\left\lvert Z_j\right\rvert\left\lvert Z_k\right\rvert\right)
\leq\sigma^2\left(\sum_{j=1}^\infty\left\lVert B^{-j}\right\rVert_F\right)^2\\
&\leq\sigma^2C\left(B^{-1}\right)^2\left(\sum_{j=1}^\infty\left\lVert B^{-j}\right\rVert_{B^{-1}}\right)^2
\leq\sigma^2C\left(B^{-1}\right)^2\left(\sum_{j=1}^\infty\left\lVert B^{-1}\right\rVert_{B^{-1}}^j\right)^2\\
&<\infty\,.
\end{align*}

\begin{Lemma}
\label{Lemma:Stationarity}
The following processes are {\rm(}strictly{\rm)} stationary and ergodic:\smallskip

\noindent\parbox{25pt}{\rm(i)} $\left(U_n\right)_{n\geq0}$\,,\smallskip

\noindent\parbox{25pt}{\rm(ii)} $\left(Y_{n-d+1}\right)_{n\geq0}$\,,\smallskip

\noindent\parbox{25pt}{\rm(iii)} $\left(U_{n-1}Y_n\right)_{n\geq1}$\,,\smallskip

\noindent\parbox{25pt}{\rm(iv)} $\left(h\left(U_n\right)Z_n\right)_{n\geq1}$ for some Borel measurable map
$h:\RR^d\rightarrow\RR^m$ {\rm(}$m\in\NN${\rm)}\,.\smallskip
\end{Lemma}
\medskip

\textit{Proof.} Let $Z=\left(Z_j\right)_{j\geq1}$ and let $S:\RR^\NN\rightarrow\RR^\NN$ denote the
shift defined by $S\left(x\right)=S\left(\left(x_j\right)_{j\geq1}\right)=\left(x_{j+1}\right)_{j\geq1}$.
Since $Z$ is stationary and ergodic, it is enough to establish an almost sure representation of the form
$f\left(S^n\left(Z\right)\right), n\in T$, if $T\in\left\{\NN,\NN_0\right\}$ or
$f\left(S^{n-1}\left(Z\right)\right), n\in T$, if $T=\NN$ for the above processes, where 
$f:\RR^\NN\rightarrow\RR^d$ is a Borel measurable map (see e.g. Lemmas~A.1 and A.2 in the Appendix of
\cite{HaeuslerLuschgy2}).
\medskip

(i) Define
\begin{displaymath}
\Lambda=\left\{x\in\RR^\NN:\sum_{j=1}^\infty\left\lVert B^{-1}e_1\right\rVert\left\lvert x_j\right\rvert<\infty
\right\}\in\B\left(\RR^\NN\right)
\end{displaymath}
and $f_1:\RR^\NN\rightarrow\RR^d$ by
\begin{displaymath}
f_1\left(x\right)=\left\{
\begin{array}{lcl}
-\sum_{j=1}^\infty B^{-j}e_1x_j & , & x\in\Lambda\\
0                               & , & x\notin\Lambda\,.
\end{array}\right.
\end{displaymath}
Then $f_1$ is $\left(\B\left(\RR^\NN\right),\B\left(\RR^d\right)\right)$-measurable, and for $n\geq0$ we have
\begin{displaymath}
U_n=-\sum_{k=1}^\infty B^{-k}e_1Z_{n+k}=f_1\left(S^n\left(Z\right)\right)\quad\text{almost surely}\,.
\end{displaymath}
\medskip

(ii) Let $f_2=\pi_d\circ f_1$. Then $f_2\left(S^n\left(Z\right)\right)=\pi_d\left(U_n\right)=Y_{n-d+1}$ almost surely,
$n\geq0$.
\medskip

(iii) Let $f_3\left(x\right)=f_1\left(x\right)\cdot\pi_1\circ f_1\left(S\left(x\right)\right)$. Then
\begin{displaymath}
f_3\left(S^{n-1}\left(Z\right)\right)=f_1\left(S^{n-1}\left(Z\right)\right)\cdot\pi_1\circ f_1\left(S^n\left(Z\right)\right)
=U_{n-1}Y_n\quad\text{almost surely},\; n\geq1\,.
\end{displaymath}

(iv) Let $f_4\left(x\right)=x_1\cdot h\circ f_1\left(S\left(x\right)\right)$. Then for $n\geq1$, almost surely,
\begin{align*}
f_4\left(S^{n-1}\left(Z\right)\right)&=S^{n-1}\left(Z\right)_1\cdot h\circ f_1\left(S\left(S^{n-1}\left(Z\right)\right)\right)\\
&=\left(\left(Z_{n-1+j}\right)_{j\geq1}\right)_1\cdot h\circ f_1\left(S^n\left(Z\right)\right)
=Z_nh\left(U_n\right)\,.
\end{align*}

\vspace{-37pt}

\phantom{x}\hfill$\Box$

\phantom{x}

\smallskip

The almost sure uniqueness of the stationary solution of \eqref{Eq:Unautoregression} being already more or less known is formulated
explicitly in the following lemma.
\medskip

\begin{Lemma}
\label{Lemma:Uniqueness}
Let $\left(\overline{U}_n\right)_{n\geq0}$ be a solution of \eqref{Eq:Unautoregression} consisting of $\RR^d$-valued random vectors
with $\left\lVert\overline{U}_0\right\rVert\in\LL^1\left(P\right)$. The following assertions are equivalent:
\smallskip

\noindent\parbox{25pt}{\rm(i)} $\left(\overline{U}_n\right)_{n\geq0}$ is stationary\,,\smallskip

\noindent\parbox{25pt}{\rm(ii)} $B^{-n}\overline{U}_n\rightarrow0$\quad almost surely as $n\to\infty$\,,\smallskip

\noindent\parbox{25pt}{\rm(iii)} $\overline{U}_0=-\sum_{j=1}^\infty B^{-j}W_j$\quad almost surely.

\end{Lemma}
\medskip

\textit{Proof.} (i)\,$\Rightarrow$\,(ii). By stationarity, we have
\begin{align*}
\sum_{n=0}^\infty E\left(\left\lVert B^{-n}\overline{U}_n\right\rVert\right)&\leq\sum_{n=0}^\infty\left\lVert B^{-n}\right\rVert_F
E\left(\left\lVert\overline{U}_0\right\rVert\right)\leq C\left(B^{-1}\right) E\left(\left\lVert\overline{U}_0\right\rVert\right)
\sum_{n=0}^\infty\left\lVert B^{-n}\right\rVert_{B^{-1}}\\
&\leq C\left(B^{-1}\right) E\left(\left\lVert\overline{U}_0\right\rVert\right)
\sum_{n=0}^\infty\left\lVert B^{-1}\right\rVert_{B^{-1}}^n<\infty
\end{align*}
which yields (ii).
\smallskip

(ii)\,$\Rightarrow$\,(iii). We know from \eqref{Eq:B-nUnconvergence} that
\begin{displaymath}
B^{-n}\overline{U}_n=\overline{U}_0+\sum_{j=1}^nB^{-j}W_j\rightarrow\overline{U}_0+\sum_{j=1}^\infty B^{-j}W_j\quad
\text{almost surely as }n\to\infty
\end{displaymath}
which implies (iii).
\smallskip

(iii)\,$\Rightarrow$\,(i) follows from Lemma~\ref{Lemma:Stationarity}~(i).\hfill$\Box$
\medskip

\begin{Remark}\label{Remark:boundedsolutions}
Any solution $\left(\overline{U}_n\right)_{n\geq0}$ of \eqref{Eq:Unautoregression} consisting of
$\RR^d$-valued random vectors (even without $\left\lVert U_0\right\rVert\in\LL^1\left(P\right)$) satisfies
\begin{displaymath}
B^{-n}\overline{U}_n\rightarrow\overline{U}_0+\sum_{j=1}^\infty B^{-j}W_j\quad\text{almost surely as }n\to\infty
\end{displaymath}
by \eqref{Eq:B-nUnconvergence}. If $\left(\overline{U}_n\right)_{n\geq0}$ is bounded in probability, then
\begin{displaymath}
B^{-n}\overline{U}_n\rightarrow0\quad\text{in probability as }n\to\infty
\end{displaymath}
because of $\rho\left(B^{-1}\right)<\infty$, and we obtain $\overline{U}_0=-\sum_{j=1}^\infty B^{-j}W_j$ almost surely,
which means that $\left(\overline{U}_n\right)_{n\geq0}$ equals almost surely the unique stationary solution of
\eqref{Eq:Unautoregression} given by \eqref{Eq:RepresentationUn}. Thus, the stationary solution of 
\eqref{Eq:Unautoregression} is the only solution which is bounded in probability.\hfill$\Box$
\end{Remark} 
\medskip
 
From now on $\left(U_n\right)_{n\geq0}$ always denotes the stationary and ergodic process given
by~\eqref{Eq:RepresentationUn}. We have 
\begin{displaymath}
\frac{1}{n}\sum_{k=1}^nU_k\rightarrow E\left(U_1\right)=0\quad\text{almost surely as }n\to\infty
\end{displaymath}
by the pointwise ergodic theorem. Now we establish a central limit theorem for $U_n$ 
though $\left(U_n\right)_{n\geq1}$ is not a martingale difference sequence.
Note that
\begin{equation}
\label{Eq:Bhochminuseins}
B^{-1}=B^{-1}\left(\theta\right)=
\begin{pmatrix} 
0                  &                            &         &                             \\
\vdots             &                            & I_{d-1} &                             \\ 
0                  &                            &         &                             \\
\frac{1}{\theta_d} & -\frac{\theta_1}{\theta_d} & \ldots  & -\frac{\theta_{d-1}}{\theta_d}        
\end{pmatrix}.
\end{equation}
\medskip

\begin{Proposition}
\label{Proposition:CLTforUn}
{\rm(}Mixing CLT for $U_n$ and $Y_n${\rm)} We have
\begin{displaymath}
\frac{1}{\sqrt{n}}\sum_{k=1}^nU_k\rightarrow\frac{\sigma}{\sum_{i=1}^d\theta_i-1}
\begin{pmatrix}
1     \\
1     \\
\vdots\\
1
\end{pmatrix}
N_1\quad\G\text{-mixing as }n\to\infty\,,
\end{displaymath}
where $\G=\sigma\left(Z_j,j\geq1\right)$, $P^{N_1}=N\left(0,1\right)$ and $N_1$ is independent of $\G$.
In particular,
\begin{displaymath}
\frac{1}{\sqrt{n}}\sum_{k=1}^nY_k\rightarrow\frac{\sigma}{\sum_{i=1}^d\theta_i-1}
N_1\quad\G\text{-mixing as }n\to\infty\,.
\end{displaymath} 
\end{Proposition}
\medskip

The above statements may be read as 
\begin{displaymath}
\frac{1}{\sqrt{n}}\sum_{k=1}^nU_k\rightarrow N\left(0,\frac{\sigma^2}{\left(\sum_{i=1}^d\theta_i-1\right)^2}E_d\right)
\quad\text{mixing}
\end{displaymath}
and
\begin{displaymath}
\frac{1}{\sqrt{n}}\sum_{k=1}^nY_k\rightarrow N\left(0,\frac{\sigma^2}{\left(\sum_{i=1}^d\theta_i-1\right)^2}\right)
\quad\text{mixing}\,,
\end{displaymath}
respectively, where $E_d$ denotes the $d\times d$-matrix with all entries equal to one and mixing means $\F$-mixing.
For the notion of mixing (and more generally, stable) convergence we refer the reader to \cite{HaeuslerLuschgy1}.
\medskip

\textit{Proof.} From \eqref{Eq:Unautoregression} we obtain $B^{-1}U_k=U_{k-1}+B^{-1}W_k,\;k\geq1$\,, so that for every $n\geq1$
\begin{displaymath}
B^{-1}\sum_{k=1}^nU_k=\sum_{k=1}^nU_{k-1}+B^{-1}\sum_{k=1}^nW_k=\sum_{k=1}^nU_k-U_n+U_0+B^{-1}e_1\sum_{k=1}^nZ_k\,.
\end{displaymath}
Consequently, using $B^{-1}e_1=\frac{1}{\theta_d}e_d$\,,
\begin{displaymath}
\left(I_d-B^{-1}\right)\sum_{k=1}^nU_k=U_n-U_0-\frac{1}{\theta_d}e_d\sum_{k=1}^nZ_k
\end{displaymath}
and hence
\begin{displaymath}
\frac{1}{\sqrt{n}}\sum_{k=1}^nU_k=\left(I_d-B^{-1}\right)^{-1}\frac{1}{\sqrt{n}}\left(U_n-U_0\right)-
\frac{1}{\theta_d}\left(I_d-B^{-1}\right)^{-1}e_d\frac{1}{\sqrt{n}}\sum_{k=1}^nZ_k\,.
\end{displaymath}
(Note that the matrix $I_d-B^{-1}$ is invertible with $\left(I_d-B^{-1}\right)^{-1}=
\sum_{j=0}^\infty B^{-j}$.) Now it is classically known that
\begin{displaymath}
-\frac{1}{\sqrt{n}}\sum_{k=1}^nZ_k\rightarrow\sigma N_1\quad\G\text{-mixing as }n\to\infty
\end{displaymath}
(see \cite{HaeuslerLuschgy1}, Example~3.16). Moreover, we have
\begin{displaymath}
\frac{1}{\sqrt{n}}\left(U_n-U_0\right)\rightarrow0\quad\text{in probability as }n\to\infty
\end{displaymath}
by stationarity of $\left(U_n\right)_{n\geq0}$ (see Lemma~\ref{Lemma:Stationarity}~(i)). This yields
\begin{displaymath}
\frac{1}{\sqrt{n}}\sum_{k=1}^nU_k\rightarrow\frac{\sigma}{\theta_d}\left(I_d-B^{-1}\right)^{-1}e_dN_1\quad
\G\text{-mixing as }n\to\infty
\end{displaymath}
(see~\cite{HaeuslerLuschgy1}, Theorem~3.18~(a),~(c)). Since
\begin{displaymath}
I_d-B^{-1}=
\begin{pmatrix}
1                   & -1                        & 0      &  \multicolumn{4}{c}{\ \ \ \ldots}                        & 0      \\
0                   &  1                        & -1     & 0 & \multicolumn{3}{c}{\quad\ \ldots}                     & 0      \\
\vdots              &                           &        &   & &              &                               & \vdots \\ 
0                   &                           & \ldots &   & & 0            & 1                             & -1     \\
-\frac{1}{\theta_d} & \frac{\theta_1}{\theta_d} & \multicolumn{4}{c}{\dots}   & \frac{\theta_{d-2}}{\theta_d} & 1+\frac{\theta_{d-1}}{\theta_d}
\end{pmatrix}
\end{displaymath}   
one easily checks that the $d$-th column of $\left(I_d-B^{-1}\right)^{-1}$ is given by
\begin{displaymath}
\left(I_d-B^{-1}\right)^{-1}e_d=\left(\left(I_d-B^{-1}\right)^{-1}\right)_{{\scriptscriptstyle\bullet}d}=
\frac{\theta_d}{\sum_{i=1}^d\theta_i-1}
\begin{pmatrix}
1 \\
1\\
\vdots \\
1
\end{pmatrix}
\end{displaymath}
(where invertibility of $I_d-B^{-1}$ assures that $\sum_{i=1}^d\theta_i-1\neq0$; for $d=2$ see also Remark~\ref{Remark:d=2})
so that the proof of the assertion concerning $U_n$ is complete. Since $Y_k=\pi_1\left(U_k\right)$, 
the assertion for $Y_n$ follows from \cite{HaeuslerLuschgy1},
Theorem~3.18 (c).\hfill$\Box$
\medskip

The covariance matrix of $U_0$ is given by $\Gamma\defeq E\left(U_0U_0^T\right)=\sigma^2\Sigma$, where
\begin{equation}
\label{Eq:Sigma}
\Sigma=\sum_{j=1}^\infty B^{-j}\widetilde{I}_d\left(B^{-j}\right)^T
\end{equation}
and $\widetilde{I}_d=e_1e_1^T$. Clearly, $\Gamma$ and $\Sigma$ are symmetric and positive semidefinite.
\medskip 

\begin{Lemma}
\label{Lemma:Sigmapositivedefinite}
The matrix $\sum_{j=1}^d B^{-j}\widetilde{I}_d\left(B^{-j}\right)^T$ is symmetric and positive definite.
Consequently, $\Sigma$ is positive definite and hence invertible.
\end{Lemma}
\medskip

\textit{Proof.} Clearly, $\sum_{j=1}^d B^{-j}\widetilde{I}_d\left(B^{-j}\right)^T$ and $\Sigma$  are symmetric
and positive semidefinite. Since $\widetilde{I}_d$ is symmetric and idempotent, we have for every $j\geq1$
and $u\in\RR^d$
\begin{displaymath}
u^TB^{-j}\widetilde{I}_d\left(B^{-j}\right)^Tu=\left\lVert\left(B^{-j}\widetilde{I}_d\right)^Tu\right\rVert^2
=\left\langle\left(B^{-j}\right)_{\scriptscriptstyle\bullet1},u\right\rangle^2\,,
\end{displaymath}
where $\left(B^{-j}\right)_{\scriptscriptstyle\bullet1}$ denotes the first column of $B^{-j}$. Consequently,
\begin{displaymath}
u^T\left(\sum_{j=1}^d B^{-j}\widetilde{I}_d\left(B^{-j}\right)^T\right)u
=\sum_{j=1}^d u^TB^{-j}\widetilde{I}_d\left(B^{-j}\right)^Tu
=\sum_{j=1}^d\left\langle\left(B^{-j}\right)_{\scriptscriptstyle\bullet1},u\right\rangle^2\,.
\end{displaymath}
For every $j=1,\ldots,d$, the first column of $B^{-j}$ satisfies
\begin{equation}
\label{Eq:Firstcolumn}
\left(B^{-j}\right)_{k1}=0\qquad\text{and}\qquad\left(B^{-j}\right)_{d-j+1,1}=\frac{1}{\theta_d}
\end{equation}
for every $k=1,\ldots,d-j$. This follows by induction from \eqref{Eq:Bhochminuseins}
so that $\left(B^{-1}\right)_{\scriptscriptstyle\bullet1}=\left(0,\ldots,0,1/\theta_d\right)^T$ and
\begin{align*}
\left(B^{-\left(j+1\right)}\right)_{\scriptscriptstyle\bullet1}&=\left(B^{-1}B^{-j}\right)_{\scriptscriptstyle\bullet1}
=\left(\left\langle\left(B^{-1}\right)_{1\scriptscriptstyle\bullet},\left(B^{-j}\right)_{\scriptscriptstyle\bullet1}\right\rangle,
\ldots,\left\langle\left(B^{-1}\right)_{d\scriptscriptstyle\bullet},\left(B^{-j}\right)_{\scriptscriptstyle\bullet1}\right\rangle\right)^T\\
&=\left(\left(B^{-j}\right)_{21},\ldots,\left(B^{-j}\right)_{d1},\frac{1}{\theta_d}\left(B^{-j}\right)_{11}-
\sum_{m=1}^{d-1}\frac{\theta_m}{\theta_d}\left(B^{-j}\right)_{m+1,1}\right)^T\,,
\end{align*}
where $\left(B^{-1}\right)_{m\scriptscriptstyle\bullet}$ denotes the $m$-th row of $B^{-1}$.
\smallskip

Now assume $\sum_{j=1}^d\left\langle\left(B^{-j}\right)_{\scriptscriptstyle\bullet1},u\right\rangle^2=0$ for some
$u\in\RR^d$. Hence $\left\langle\left(B^{-j}\right)_{\scriptscriptstyle\bullet1},u\right\rangle=0$ for every
$j=1,\ldots,d$ or equivalently $Du=0$, where the $d\times d$-matrix $D$ is given by
\begin{displaymath}
D=
\begin{pmatrix}
\left(B^{-d}\right)_{\scriptscriptstyle\bullet1}^T \\
\vdots                                             \\
\left(B^{-1}\right)_{\scriptscriptstyle\bullet1}^T
\end{pmatrix}\,.
\end{displaymath}
Then $D$ is an upper triangular matrix with all diagonal elements equal to $1/\theta_d$ so that
$\det\left(D\right)=\left(1/\theta_d\right)^d\neq0$. We conclude that $u=0$ is the only solution
of the equation $Du=0$.\hfill$\Box$
\medskip

Let $\gamma\left(k\right)\defeq E\left(Y_0Y_k\right)$ for $k\geq0$. Then by stationarity of
$\left(Y_{n-d+1}\right)_{n\geq0}$ (see Lemma~\ref{Lemma:Stationarity}~(ii)), we obtain 
$E\left(Y_jY_k\right)=\gamma\left(\left\lvert j-k\right\rvert\right)$ for $j,k\geq-d+1$ so that
\begin{displaymath}
\Gamma=\left(\gamma\left(\left\lvert j-k\right\rvert\right)\right)_{1\leq j,k\leq d}=
\begin{pmatrix}
\gamma\left(0\right)   & \gamma\left(1\right)   & \ldots & \gamma\left(d-1\right) \\
\gamma\left(1\right)   & \gamma\left(0\right)   & \ldots & \gamma\left(d-2\right) \\
\vdots                 &                        &        & \vdots                 \\
\gamma\left(d-1\right) & \gamma\left(d-2\right) & \ldots & \gamma\left(0\right)   \\
\end{pmatrix}.
\end{displaymath}
In our setting, the restricted Yule-Walker equations hold true, namely
\begin{equation}
\label{Eq:RestrictedYuleWalker}
\gamma\left(k\right)=\sum_{j=1}^d\gamma\left(\left\lvert j-k\right\rvert\right)\theta_j\quad
\text{for }1\leq k\leq d-1\,.
\end{equation}
In fact, for every $k\geq1$, $U_k=\left(Y_k,Y_{k-1},\ldots,Y_{k-d+1}\right)^T$ and $Z_k$ are independent
by \eqref{Eq:RepresentationUn}, hence $Y_{k-m}$ and $Z_k$ are independent for $0\leq m\leq d-1$ so that
$E\left(Y_{k-m}Z_k\right)=E\left(Y_{k-m}\right)E\left(Z_k\right)=0$. For $1\leq k\leq d-1$, we obtain
\begin{align*}
\gamma\left(k\right)&=E\left(Y_0Y_k\right)=E\left(Y_0\left(\theta_1Y_{k-1}+\cdots+\theta_dY_{k-d}+Z_k\right)\right)\\
&=\sum_{j=1}^d\theta_jE\left(Y_0Y_{k-j}\right)+E\left(Y_0Z_k\right)
=\sum_{j=1}^d\theta_j\gamma\left(\left\lvert k-j\right\rvert\right)
\end{align*}
using $Y_0=Y_{k-k}$. For $k=d$, we have
\begin{equation}
\label{Eq:Gammad}
\gamma\left(d\right)=\sum_{j=1}^d\gamma\left(\left\lvert j-d\right\rvert\right)\theta_j-\frac{\sigma^2}{\theta_d}
\end{equation}
since
\begin{align*}
\gamma\left(d\right)&=E\left(Y_0Y_d\right)=E\left(Y_0\left(\theta_1Y_{d-1}+\cdots+\theta_dY_0+Z_d\right)\right)\\
&=\sum_{j=1}^d\theta_j\gamma\left(d-j\right)+E\left(Y_0Z_d\right)
\end{align*}
and by \eqref{Eq:Firstcolumn}, $\left(B^{-1}\right)_{11}=\cdots=\left(B^{-\left(d-1\right)}\right)_{11}=0$ and
$\left(B^{-d}\right)_{11}=1/\theta_d$ so that
\begin{displaymath}
E\left(Y_0Z_d\right)=-\sum_{j=1}^\infty\left(B^{-j}\right)_{11}E\left(Z_jZ_d\right)
=-\left(B^{-d}\right)_{11}E\left(Z_d^2\right)=-\frac{\sigma^2}{\theta_d}
\end{displaymath}
because $E\left(Z_jZ_d\right)=E\left(Z_j\right)E\left(Z_d\right)=0$ for $j\neq d$. On the other
hand, we have for $k=0$ 
\begin{align}
\gamma\left(0\right)&=E\left(Y_d^2\right)=
E\left(Y_d\left(\theta_1Y_{d-1}+\cdots+\theta_dY_0+Z_d\right)\right)\nonumber\\
&=\theta_1\gamma\left(1\right)+\cdots+\theta_d\gamma\left(d\right)\label{Eq:Gamma0}
\end{align}
so that
\begin{displaymath}
\gamma\left(d\right)=\gamma\left(0\right)\frac{1}{\theta_d}-\gamma\left(1\right)\frac{\theta_1}{\theta_d}-
\cdots-\gamma\left(d-1\right)\frac{\theta_{d-1}}{\theta_d}\,.
\end{displaymath}
This suggests to define
\begin{equation}
\label{Eq:thetastar}
\theta^\ast\defeq\left(-\frac{\theta_{d-1}}{\theta_d},\ldots,-\frac{\theta_1}{\theta_d},\frac{1}{\theta_d}\right)^T
\end{equation}
which yields
\begin{equation}
\label{Eq:YuleWalkerd}
\gamma\left(d\right)=\sum_{j=1}^d\gamma\left(\left\lvert j-d\right\rvert\right)\theta_j^\ast\,.
\end{equation}
Now it is easy to see that the full Yule-Walker equations hold true when $\theta$ is replaced by $\theta^\ast$.
\medskip

\begin{Lemma}{\rm(}Yule-Walker equations{\rm)}
\label{Lemma:YuleWalkerequationsthetastar}
We have
\begin{displaymath}
\Gamma\theta^\ast=\left(\gamma\left(1\right),\ldots,\gamma\left(d\right)\right)^T\,.
\end{displaymath}
\end{Lemma}
\medskip

\textit{Proof.} The equations \eqref{Eq:RestrictedYuleWalker} imply for $1\leq k\leq d-1$,
\begin{displaymath}
\gamma\left(d-k\right)\theta_d=\gamma\left(k\right)-\sum_{j=1}^{d-1}\gamma\left(\left\lvert j-k\right\rvert\right)\theta_j
\end{displaymath}
or equivalently
\begin{displaymath}
\gamma\left(d-k\right)=\gamma\left(k\right)\frac{1}{\theta_d}-\sum_{j=1}^{d-1}\gamma\left(\left\lvert j-k\right\rvert\right)
\frac{\theta_j}{\theta_d}
=\gamma\left(k\right)\theta_d^\ast+\sum_{j=1}^{d-1}\gamma\left(\left\lvert j-k\right\rvert\right)\theta_{d-j}^\ast\,.
\end{displaymath}
Setting $k=d-m$ with $1\leq m\leq d-1$ gives 
\begin{displaymath}
\gamma\left(m\right)
=\gamma\left(d-m\right)\theta_d^\ast+\sum_{j=1}^{d-1}\gamma\left(\left\lvert j-d+m\right\rvert\right)\theta_{d-j}^\ast
\end{displaymath}
and finally, setting $j=d-r$ with $1\leq r\leq d-1$ yields
\begin{displaymath}
\gamma\left(m\right)
=\gamma\left(d-m\right)\theta_d^\ast+\sum_{r=1}^{d-1}\gamma\left(\left\lvert m-r\right\rvert\right)\theta_r^\ast
=\sum_{s=1}^d\gamma\left(\left\lvert s-m\right\rvert\right)\theta_s^\ast\,.
\end{displaymath} 
For the case $m=d$ see \eqref{Eq:YuleWalkerd}.\hfill$\Box$
\medskip

An immediate consequence is as follows. Since by \eqref{Eq:RestrictedYuleWalker} and \eqref{Eq:Gammad},
\begin{displaymath}
\Gamma\theta=\left(\gamma\left(1\right),\ldots,\gamma\left(d-1\right),\gamma\left(d\right)+\frac{\sigma^2}{\theta_d}\right)^T\,,
\end{displaymath}
one obtains $\Gamma\theta^\ast-\Gamma\theta=-\left(\sigma^2/\theta_d\right)e_d$ or
\begin{equation}
\label{Eq:thetastarminustheta}
\theta^\ast-\theta=-\frac{\sigma^2}{\theta_d}\Gamma^{-1}e_d=-\frac{1}{\theta_d}\Sigma^{-1}e_d\,.
\end{equation}
\medskip

Next we observe some strong laws which are consequences of the Birkhoff ergodic theorem.
\medskip

\begin{Lemma}
\label{Lemma:Stronglaws}
Let $h:\RR^d\rightarrow\RR^m$ be a Borel measurable map. As $n\to\infty$, we have

\noindent\parbox{25pt}{\rm(a)}\quad $\displaystyle\frac{1}{n}\sum_{k=0}^{n-1}h\left(U_k\right)h\left(U_k\right)^T
\rightarrow E\left(h\left(U_0\right)h\left(U_0\right)^T\right)$\quad almost surely,

provided $\left\lVert h\left(U_0\right)\right\rVert\in\LL^2\left(P\right)$ and, in particular,
\begin{displaymath}
\frac{1}{n}\sum_{k=0}^{n-1}U_kU_k^T\rightarrow\Gamma\quad\text{almost surely}\,,
\end{displaymath}

\noindent\parbox{25pt}{\rm(b)}\quad $\displaystyle\frac{1}{n}\sum_{k=1}^nU_{k-1}Y_k\rightarrow\Gamma\theta^\ast$\quad almost surely,

\noindent\parbox{25pt}{\rm(c)}\quad $\displaystyle\frac{1}{n}\sum_{k=1}^nh\left(U_k\right)Z_k\rightarrow0$\quad almost surely

provided $\left\lVert h\left(U_1\right)\right\rVert\in\LL^1\left(P\right)$.
\end{Lemma}
\medskip

\textit{Proof.} Using $E\left(\left\lVert h\left(U_1\right)h\left(U_1\right)^T\right\rVert_F\right)=
E\left(\left\lVert h\left(U_1\right)\right\rVert^2\right)<\infty$, assertion (a)
follows from Lemma~\ref{Lemma:Stationarity}~(i) and the ergodic theorem. Likewise,
it follows from Lemma~\ref{Lemma:Stationarity}~(iii) and the ergodic theorem that, as $n\to\infty$,
\begin{displaymath}
\frac{1}{n}\sum_{k=1}^nU_{k-1}Y_k\rightarrow E\left(U_0Y_1\right)\quad\text{almost surely}
\end{displaymath}
and $E\left(U_0Y_1\right)=\left(\gamma\left(1\right),\ldots,\gamma\left(d\right)\right)^T=\Gamma\theta^\ast$
by Lemma~\ref{Lemma:YuleWalkerequationsthetastar}, and by Lemma~\ref{Lemma:Stationarity}~(iv)
\begin{displaymath}
\frac{1}{n}\sum_{k=1}^nh\left(U_k\right)Z_k\rightarrow E\left(h\left(U_1\right)Z_1\right)=0\quad\text{almost surely}
\end{displaymath}
because $U_1$ and $Z_1$ are independent. Thus, (a) - (c) have been proven.
\hfill$\Box$
\medskip

\section{\Large Central limit theorems via reverse martingale difference sequences}
\label{Section:CLT}

In the setting of the present paper a central limit theorem for reverse martingale difference sequences
turns out to be a main tool to obtain central limit theorems for processes like
$\left(U_nZ_n\right)_{n\geq0}$. For this, we consider a sequence $\left(\HH_k\right)_{k\geq1}$ of sub-$\sigma$-fields
of $\F$ with $\HH_k\supset\HH_{k+1}$ for every $k\geq1$, which we call a \textit{reverse filtration}.
Set $\HH_\infty:=\bigcap_{k\geq1}\HH_k$. We call a sequence $\left(X_k\right)_{k\geq1}$ of integrable 
$d$-dimensional random vectors which is adapted to $\left(\HH_k\right)_{k\geq1}$ (i.e. $X_k$ is $\HH_k$-measurable
for all $k\geq1$) a \textit{reverse martingale difference sequence} if $E\left(X_k|\HH_{k+1}\right)=0$ for all
$k\geq1$. The following limit theorem for the partial sums $\sum_{k=1}^nX_k$ should not be confused with central 
limit theorems for reverse martingales, that is for series remainders $\sum_{k=n}^\infty X_k$ as $n\to\infty$ provided
the series $\sum_{k=1}^\infty X_k$ converges in $\LL^1\left(P\right)$
as obtained e.g. in \cite{Loynes} for $d=1$. It is a version of Corollary~3.4 in \cite{HaeuslerLuschgy2} using a
special real norming sequence.
\medskip

\begin{Theorem}
\label{Theorem:Reversedmartingaledifferences}
Let $\left(X_k\right)_{k\geq1}$ be a square-inte\-grable $\RR^d$-valued reverse martingale difference sequence w.r.t. 
the reverse filtration $\left(\HH_k\right)_{k\geq1}$. If
\begin{align}
&\frac{1}{n}\sum_{k=1}^nE\left(\left\lVert X_k\right\rVert^21_{\{\left\lVert X_k\right\rVert\geq\varepsilon\sqrt{n}\}}|
\HH_{k+1}\right)\rightarrow0\quad\text{in probability as }n\to\infty\label{Condition:CLB}\\
&\text{for every }\varepsilon>0\nonumber
\end{align}
and
\begin{align}
&\frac{1}{n}\sum_{k=1}^nE\left(X_kX_k^T|\HH_{k+1}\right)\rightarrow A\quad\text{in probability as }
n\to\infty\text{ for some}\label{Condition:norming}\\
&\HH_\infty\text{-measurable symmetric and positive semi-definite random}\nonumber\\
&d\times d\text{-matrix }A\,,\nonumber
\end{align}
then
\begin{displaymath}
\frac{1}{\sqrt{n}}\sum_{k=1}^nX_k\rightarrow A^{1/2}N_d\quad\HH_\infty\text{-stably as }n\to\infty\,,
\end{displaymath}
where $P^{N_d}=N\left(0,I_d\right)$ and $N_d$ is independent of $\HH_\infty$.
\end{Theorem}

The convergence in Theorem~\ref{Theorem:Reversedmartingaledifferences} is mixing if the matrix $A$ is nonrandom.
\medskip

In the setting of the previous sections, $\left(\G_n\right)_{n\geq1}$ with $\G_n\defeq\sigma\left(Z_j:j\geq n\right)$
for every $n\geq1$ is a reverse filtration. Let $h:\RR^d\rightarrow\RR^m$ for some $m\geq1$ be Borel measurable. Then
$h\left(U_n\right)Z_n$ is $\G_n$-measurable for every $n\geq1$. If $\left\lVert h\left(U_0\right)\right\rVert\in\LL^2\left(P\right)$,
then, by independence of $U_n$ and $Z_n$, for the $m$-dimensional random vector $h\left(U_n\right)Z_n$ we have
$\left\lVert h\left(U_n\right)Z_n\right\rVert\in\LL^2\left(P\right)$ and, by $\G_{n+1}$-measurability of $U_n$
and independence of $Z_n$ and $\G_{n+1}$,
\begin{displaymath}
E\left(h\left(U_n\right)Z_n|\G_{n+1}\right)=h\left(U_n\right)E\left(Z_n|\G_{n+1}\right)=h\left(U_n\right)E\left(Z_n\right)=0\,.
\end{displaymath}
Consequently, $\left(h\left(U_n\right)Z_n\right)_{n\geq1}$ is a square-integrable reverse martingale difference sequence
w.r.t. $\left(\G_n\right)_{n\geq1}$, and Theorem~\ref{Theorem:Reversedmartingaledifferences} implies the following
result. However, since $\G_\infty=\bigcap_{n=1}^\infty\G_n$ is $P$-trivial by the Kolmogorov $0-1$ law,   
$\G_\infty$-stable convergence reduces to mere distributional convergence (cf.\cite{HaeuslerLuschgy1}, Exercise~3.4).
\medskip

\begin{Proposition} 
\label{Proposition:Clt}
If $\left\lVert h\left(U_0\right)\right\rVert\in\LL^2\left(P\right)$, then
\begin{displaymath}
\frac{1}{\sqrt{n}}\sum_{k=1}^nh\left(U_k\right)Z_k\stackrel{\D}{\rightarrow}
N\left(0,\sigma^2E\left(h\left(U_1\right)h\left(U_1\right)^T\right)\right)
\end{displaymath}
{\rm(}where $\stackrel{\D}{\rightarrow}$ means convergence in distribution{\rm)}.
\end{Proposition}
\medskip

\textit{Proof.} We have to verify~\eqref{Condition:CLB} and~\eqref{Condition:norming} for $X_k=h\left(U_k\right)Z_k$.
By stationarity of $\left(h\left(U_n\right)Z_n\right)_{n\geq1}$, see Lemma~\ref{Lemma:Stationarity}~(iv), for every $\varepsilon>0$
\begin{align*}
&E\left(\frac{1}{n}\sum_{k=1}^nE\left(\left\lVert h\left(U_k\right)Z_k\right\rVert^2
1_{\{\left\lVert h\left(U_k\right)Z_k\right\rVert\geq\varepsilon\sqrt{n}\}}|\G_{k+1}\right)\right)\\
&\qquad=\frac{1}{n}\sum_{k=1}^nE\left(\left\lVert h\left(U_k\right)Z_k\right\rVert^2
1_{\{\left\lVert h\left(U_k\right)Z_k\right\rVert\geq\varepsilon\sqrt{n}\}}\right)\\
&\qquad=E\left(\left\lVert h\left(U_1\right)Z_1\right\rVert^2
1_{\{\left\lVert h\left(U_1\right)Z_1\right\rVert\geq\varepsilon\sqrt{n}\}}\right)\rightarrow0\quad\text{as }n\to\infty
\end{align*}
since $E\left(\left\lVert h\left(U_1\right)Z_1\right\rVert^2\right)<\infty$ so that~\eqref{Condition:CLB} holds true and
\begin{align*}
&\frac{1}{n}\sum_{k=1}^nE\left(h\left(U_k\right)Z_k\left(h\left(U_k\right)Z_k\right)^T|\G_{k+1}\right)= 
\frac{1}{n}\sum_{k=1}^nh\left(U_k\right)h\left(U_k\right)^TE\left(Z_k^2|\G_{k+1}\right)\\
&=\frac{\sigma^2}{n}\sum_{k=1}^nh\left(U_k\right)h\left(U_k\right)^T\rightarrow
\sigma^2E\left(h\left(U_1\right)h\left(U_1\right)^T\right)\quad\text{almost surely as }n\to\infty
\end{align*}
by Lemma~\ref{Lemma:Stronglaws}~(a) so
that~\eqref{Condition:norming} is satisfied with $A=\sigma^2E\left(h\left(U_1\right)h\left(U_1\right)^T\right)$.
Therefore, by Theorem~\ref{Theorem:Reversedmartingaledifferences},
\begin{displaymath}
\frac{1}{\sqrt{n}}\sum_{k=1}^nh\left(U_k\right)Z_k\rightarrow\sigma E\left(h\left(U_1\right)h\left(U_1\right)^T\right)^{1/2} 
N_d\quad\G_\infty\text{-mixing as }n\to\infty\,,
\end{displaymath}
where $P^{N_d}=N\left(0,I_d\right)$ and $N_d$ is independent of $\G_\infty$. This is the
desired result.\hfill$\Box$
\medskip

As an immediate consequence of Proposition~\ref{Proposition:Clt} we obtain
\begin{displaymath}
\frac{1}{\sqrt{n}}\sum_{k=1}^nU_kZ_k\stackrel{\D}{\rightarrow}N\left(0,\sigma^2\Gamma\right)
=N\left(0,\sigma^4\Sigma\right)\quad\text{as }n\to\infty
\end{displaymath}
and in particular for $0\leq m\leq d-1$
\begin{displaymath}
\frac{1}{\sqrt{n}}\sum_{k=1}^nY_{k-m}Z_k\stackrel{\D}{\rightarrow}N\left(0,\sigma^2\gamma\left(0\right)\right)\quad
\text{as }n\to\infty\,.
\end{displaymath}
\medskip

In the remaining part of this section we consider the least squares estimator of the true parameter given by
\begin{equation}
\label{Eq:leastsquares}
\widehat{\theta}_n=\left(\sum_{k=1}^nU_{k-1}U_{k-1}^T\right)^{-1}\sum_{k=1}^nU_{k-1}Y_k\,,\quad n\geq1\,,
\end{equation}
where 
$\left(\sum_{k=1}^nU_{k-1}U_{k-1}^T\right)^{-1}:=0$ if the random symmetric and positive semidefinite matrix
$\sum_{k=1}^nU_{k-1}U_{k-1}^T$ is not positive definite and therefore singular. In our setting $\widehat{\theta}_n$
is not consistent but converges to $\theta^\ast$ and a central limit theorem holds true, as shown by the following proposition.
\medskip

\begin{Proposition}
\label{Proposition:thetanhat}
We have, as $n\to\infty$,
\begin{equation}
\label{Eq:thetanhatalmostsurely}
\widehat{\theta}_n\rightarrow\theta^\ast\quad\text{almost surely}
\end{equation}
and
\begin{equation}
\label{Eq:thetanhatindistribution}
\sqrt{n}\left(\widehat{\theta}_n-\theta^\ast\right)\stackrel{\D}{\rightarrow}
\frac{\sigma}{\theta_d}\Gamma^{-1/2}N_d\,,
\end{equation}
where $P^{N_d}=N\left(0,I_d\right)$ {\rm(}and $\Gamma^{-1/2}\defeq\left(\Gamma^{1/2}\right)^{-1}${\rm)}.
\end{Proposition}
\medskip

For the proof we rely on the following lemma.
\medskip

\begin{Lemma}
\label{Lemma:sums}
Let $a,b,u,v\in\ZZ$ with $a,b,u,v\leq d$ satisfying
\begin{equation}
\label{eq:conditionsonindices}
\left\lvert a-b\right\rvert=\left\lvert u-v\right\rvert\,.
\end{equation}
Then
\begin{displaymath}
\frac{1}{\sqrt{n}}\left(\sum_{k=1}^nY_{k-a}Y_{k-b}-\sum_{k=1}^nY_{k-u}Y_{k-v}\right)\rightarrow0\quad\text{in probability as }
n\to\infty\,.
\end{displaymath}
\end{Lemma}

In the preceeding lemma one may replace $\left(1/\sqrt{n}\right)_{n\geq1}$ by any sequence $\left(a_n\right)_{n\geq1}$
in $\RR$ with $a_n\rightarrow0$. Note that the restriction $a\leq d$ assures $k-a\geq-d+1$ so that $Y_{k-a}$ is well defined for
every $k\geq1$. The same holds true for $b,u,$ and $v$.
\medskip

\textit{Proof.} We may and will assume $a\geq b$ and $u\geq v$ so that $a-b=u-v$. Set $m-a\defeq k-u$, i.e.
$m=k+a-u$. Then by \eqref{eq:conditionsonindices} $k-v=k+a-u-b=m-b$ and for $1\leq k\leq n$, we have
$1+a-u\leq m\leq n+a-u$ so that for every $n\geq1$
\begin{equation}
\label{eq:equationsums}
\sum_{k=1}^nY_{k-u}Y_{k-v}=\sum_{m=1+a-u}^{n+a-u}Y_{m-a}Y_{m-b}\,.
\end{equation}
If $a-u<0$, this implies for every $n\geq1$
\begin{displaymath}
\sum_{k=1}^nY_{k-a}Y_{k-b}-\sum_{k=1}^nY_{k-u}Y_{k-v}=-\sum_{k=1+a-u}^0Y_{k-a}Y_{k-b}+\sum_{k=n+a-u+1}^nY_{k-a}Y_{k-b}\,.
\end{displaymath}
Clearly,
\begin{displaymath}
\frac{1}{\sqrt{n}}\sum_{k=1+a-u}^0Y_{k-a}Y_{k-b}\rightarrow0\quad\text{everywhere as }n\to\infty\,.
\end{displaymath}
By stationarity of $\left(Y_{n-d+1}\right)_{n\geq0}$ (see Lemma~\ref{Lemma:Stationarity}~(ii)), we have
\begin{align*}
&E\left(\left\lvert\frac{1}{\sqrt{n}}\sum_{k=n+a-u+1}^nY_{k-a}Y_{k-b}\right\rvert\right)\leq\frac{1}{\sqrt{n}}
\sum_{k=n+a-u+1}^nE\left(\left\lvert Y_{k-a}Y_{k-b}\right\rvert\right)\\
&\qquad\leq\frac{1}{\sqrt{n}}\sum_{k=n+a-u+1}^n\left(E\left(Y_{k-a}^2\right)\right)^{1/2}\left(E\left(Y_{k-b}^2\right)\right)^{1/2}
=\frac{1}{\sqrt{n}}\left\lvert a-u\right\rvert E\left(Y_0^2\right)\rightarrow0
\end{align*}
as $n\to\infty$ which entails
\begin{displaymath} 
\frac{1}{\sqrt{n}}\sum_{k=n+a-u+1}^nY_{k-a}Y_{k-b}\rightarrow0\quad\text{in probability as }n\to\infty\,.
\end{displaymath}
Thus the assertion.
\medskip

If $a-u=0$ then $b-v=0$ by \eqref{eq:conditionsonindices} and hence for every $n\geq1$
\begin{displaymath}
\sum_{k=1}^nY_{k-a}Y_{k-b}-\sum_{k=1}^nY_{k-u}Y_{k-v}=0\,.
\end{displaymath}
If $a-u>0$, then by \eqref{eq:equationsums}
\begin{displaymath}
\sum_{k=1}^nY_{k-a}Y_{k-b}-\sum_{k=1}^nY_{k-u}Y_{k-v}=\sum_{k=1}^{a-u}Y_{k-a}Y_{k-b}-\sum_{k=n+1}^{n+a-u}Y_{k-a}Y_{k-b}
\end{displaymath}
and the same arguments as for the first case yield the assertion. In fact, this case reduces to the first one by
interchanging $a$ and $u$ (and $b$ and $v$).\hfill$\Box$
\medskip

\textit{Proof of Proposition~\ref{Proposition:thetanhat}.} By 
Lemma~\ref{Lemma:Sigmapositivedefinite} and Lemma~\ref{Lemma:Stronglaws} (a), (b)
\begin{displaymath}
\widehat{\theta}_n\rightarrow\Gamma^{-1}\Gamma\theta^\ast=\theta^\ast\quad\text{almost surely as }n\to\infty
\end{displaymath}
so that~\eqref{Eq:thetanhatalmostsurely} holds true.
\smallskip

For the proof of~\eqref{Eq:thetanhatindistribution} we use $U_{k-1}=B^{-1}U_k-B^{-1}W_k, k\geq1$, from~\eqref{Eq:Unautoregression}
to obtain, for $n\geq1$,
\begin{align*}
&\sum_{k=1}^nU_{k-1}Y_k-\sum_{k-1}^nU_{k-1}U_{k-1}^T\theta^\ast\\
&\quad=\sum_{k=1}^n\left(B^{-1}U_k-B^{-1}W_k\right)Y_k-\sum_{k=1}^nU_{k-1}U_{k-1}^T\theta^\ast\\
&\quad=-\sum_{k=1}^nB^{-1}W_kY_k+\sum_{k=1}^n\left(B^{-1}U_kY_k-U_kU_k^T\theta^\ast\right)\\
&\qquad+\left(\sum_{k=1}^nU_kU_k^T-\sum_{k=1}^nU_{k-1}U_{k-1}^T\right)\theta^\ast\eqdef I_n+\II_n+\III_n\,.
\end{align*}
On the event $\left\{\det\left(\sum_{k=1}^nU_{k-1}U_{k-1}^T\right)>0\right\}$ we have, for $n\geq1$,
\begin{equation}
\label{Eq:representation}
\left(\frac{1}{n}\sum_{k=1}^nU_{k-1}U_{k-1}^T\right)\sqrt{n}\left(\widehat{\theta}_n-\theta^\ast\right)=
\frac{1}{\sqrt{n}}I_n+\frac{1}{\sqrt{n}}\II_n+\frac{1}{\sqrt{n}}\III_n\eqdef C_n\,,
\end{equation}
whereas all random quantities appearing in~\eqref{Eq:representation} are defined everywhere on $\Omega$. 
\smallskip

We will investigate $\frac{1}{\sqrt{n}}I_n$, $\frac{1}{\sqrt{n}}\II_n$ and $\frac{1}{\sqrt{n}}\III_n$
separately.
\smallskip

As to $\frac{1}{\sqrt{n}}\III_n$, by stationarity of $\left(U_n\right)_{n\geq0}$, see Lemma~\ref{Lemma:Stationarity}~(i),
the random matrices $\left(U_nU_n^T\right)_{n\geq0}$ are identically distributed which yields
\begin{equation}
\label{Eq:limitIIIn}
\frac{1}{\sqrt{n}}\III_n=\frac{1}{\sqrt{n}}\left(U_nU_n^T-U_0U_0^T\right)\theta^\ast\rightarrow0\quad
\text{in probability as }n\to\infty
\end{equation}
so that by the Cram\'{e}r-Slutzky theorem the summand $\frac{1}{\sqrt{n}}\III_n$ does not contribute to the
asymptotic distribution of $C_n$ in~\eqref{Eq:representation}.
\smallskip

As to $\frac{1}{\sqrt{n}}I_n$, for $n\geq1$, we have, using~\eqref{Eq:Bhochminuseins},
\begin{equation}
\label{Eq:rootIn}
\frac{1}{\sqrt{n}}I_n=-\frac{1}{\sqrt{n}}\sum_{k=1}^nB^{-1}e_1Y_kZ_k=
-\frac{1}{\theta_d}\frac{1}{\sqrt{n}}\sum_{k=1}^ne_dY_kZ_k=\frac{1}{\sqrt{n}}\sum_{k=1}^nh_I\left(U_k\right)Z_k
\end{equation}
with $h_I\left(x\right)=-\frac{1}{\theta_d}x_1e_d$, $x\in\RR^d$. 
\smallskip

For $d=1$, the situation is particularly simple. Proposition~\ref{Proposition:Clt} implies (note that 
$h_I\left(U_k\right)=h_I\left(Y_k\right)=-\frac{1}{\theta_1}Y_k$)
\begin{displaymath}
\frac{1}{\sqrt{n}}I_n\stackrel{\D}{\rightarrow}N\left(0,\frac{\sigma^2}{\theta_1^2}E\left(Y_1^2\right)\right)=
N\left(0,\frac{\sigma^2}{\theta_1^2}\gamma\left(0\right)\right)\quad\text{as }n\to\infty\,.
\end{displaymath}
Since $B^{-1}=\theta^\ast=1/\theta_1$, by~\eqref{Eq:Bhochminuseins} and~\eqref{Eq:thetastar} we have, or all $k\geq1$,
\begin{displaymath}
B^{-1}U_kY_k-U_kU_k^T\theta^\ast=\frac{1}{\theta_1}Y_k^2-Y_k^2\frac{1}{\theta_1}=0
\end{displaymath}
so that $\II_n=0$ for all $n\geq1$. Thus, in view of~\eqref{Eq:limitIIIn},
\begin{equation}
\label{Eq:limitCnd=1}
C_n\stackrel{\D}{\rightarrow}N\left(0,\frac{\sigma^2}{\theta_1^2}\gamma\left(0\right)\right)\quad\text{as }n\to\infty\,.
\end{equation}
For $d\geq2$, $\II_n$ does not disappear and we have to determine if it does contribute
to the asymptotic distribution of $C_n$ or not. Using~\eqref{Eq:Bhochminuseins}, by elementary algebra, we obtain
for $j=1,\ldots,d-1$
\begin{align}
&\left(B^{-1}U_kY_k-U_kU_k^T\theta^\ast\right)_j\nonumber\\
&\qquad=\frac{1}{\theta_d}\left(\theta_dY_kY_{k-j}+\sum_{\ell=0}^{d-2}\theta_{d-\ell-1}Y_{k-j+1}Y_{k-\ell}
-Y_{k-j+1}Y_{k-d+1}\right)\label{Eq:j}
\end{align}
and
\begin{align}
&\left(B^{-1}U_kY_k-U_kU_k^T\theta^\ast\right)_d\nonumber\\
&\qquad=\frac{1}{\theta_d}\left(Y_k^2-\sum_{\ell=1}^{d-1}\theta_\ell Y_kY_{k-\ell}+
\sum_{\ell=0}^{d-2}\theta_{d-\ell-1}Y_{k-d+1}Y_{k-\ell}
-Y_{k-d+1}^2\right)\,.\label{Eq:d}
\end{align}
We consider three cases.
\smallskip

\textit{Case 1.} $j=1$. Then from~\eqref{Eq:j}
\begin{align*}
&\theta_d\left(B^{-1}U_kY_k-U_kU_k^T\theta^\ast\right)_1\\
&\qquad=\theta_dY_kY_{k-1}+\sum_{\ell=0}^{d-2}\theta_{d-\ell-1}Y_kY_{k-\ell}-Y_kY_{k-d+1}\\
&\qquad=\theta_dY_kY_{k-1}+\sum_{\ell=0}^{d-2}\theta_{d-\ell-1}Y_kY_{k-\ell}
-\sum_{\ell=1}^d\theta_\ell Y_{k-\ell}Y_{k-d+1}-Y_{k-d+1}Z_k\\
&\qquad\eqdef A_k\left(1\right)-Y_{k-d+1}Z_k\,,
\end{align*}
where we have used~\eqref{Eq:Yautoregression}, i.e. $Y_k=\sum_{\ell=1}^d\theta_\ell Y_{k-\ell}+Z_k$. Our aim is to show
\begin{equation}
\label{Eq:A1}
\frac{1}{\sqrt{n}}\sum_{k=1}^nA_k\left(1\right)\rightarrow0\quad\text{in probability as }n\to\infty\,.
\end{equation}
For the proof of~\eqref{Eq:A1} we write
\begin{align*}
A_k\left(1\right)&=\theta_dY_kY_{k-1}+\theta_{d-1}Y_k^2+\sum_{\ell=1}^{d-2}\theta_{d-\ell-1}Y_kY_{k-\ell}
-\sum_{\ell=1}^{d-2}\theta_\ell Y_{k-\ell}Y_{k-d+1}\\
&\qquad-\theta_{d-1}Y_{k-d+1}^2-\theta_dY_{k-d}Y_{k-d+1}\\
&=\theta_d\left(Y_kY_{k-1}-Y_{k-d}Y_{k-d+1}\right)+\theta_{d-1}\left(Y_k^2-Y_{k-d+1}^2\right)\\
&\qquad+\sum_{\ell=1}^{d-2}\theta_\ell\left(Y_kY_{k-d+1+\ell}-Y_{k-\ell}Y_{k-d+1}\right)
\end{align*}
so that, for $n\geq1$,
\begin{align*}
\sum_{k=1}^nA_k\left(1\right)&=\theta_d\sum_{k=1}^n\left(Y_kY_{k-1}-Y_{k-d}Y_{k-d+1}\right)+
\theta_{d-1}\sum_{k=1}^n\left(Y_k^2-Y_{k-d+1}^2\right)\\
&\qquad+\sum_{\ell=1}^{d-2}\theta_\ell\sum_{k=1}^n\left(Y_kY_{k-d+1+\ell}-Y_{k-\ell}Y_{k-d+1}\right)\,.
\end{align*}

Now ($d$ applications of) Lemma~\ref{Lemma:sums} yields~\eqref{Eq:A1} and we get
\begin{equation}
\label{Eq:representationfirstcomponent}
\frac {1}{\sqrt{n}}\sum_{k=1}^n\left(B^{-1}U_kY_k-U_kU_k^T\theta^\ast\right)_1=-\frac{1}{\theta_d}\frac{1}{\sqrt{n}}\sum_{k=1}^n
Y_{k-d+1}Z_k+R_n\left(1\right)
\end{equation}
with $R_n\left(1\right)\rightarrow0$ in probability as $n\to\infty$.
\smallskip

\textit{Case 2.} $2\leq j\leq d-1$. Then by~\eqref{Eq:j}
\begin{align*}
&\theta_d\left(B^{-1}U_kY_k-U_kU_k^T\theta^\ast\right)_j
=\theta_dY_kY_{k-j}+\sum_{\ell=1}^{d-1}\theta_\ell Y_{k-j+1}Y_{k-d+1+\ell}-Y_{k-j+1}Y_{k-d+1}\\
&\quad=\theta_dY_kY_{k-j}+\sum_{\substack{\ell=1\\ \ell\neq j-1}}^{d-1}\theta_\ell Y_{k-j+1}Y_{k-d+1+\ell}-Y_{k-j+1}Y_{k-d+1}
+\theta_{j-1}Y_{k-j+1}Y_{k-d+j}\,.
\end{align*}
From~\eqref{Eq:Yautoregression} for $j=2,\ldots,d-1$
\begin{align*}
Y_kY_{k-d+j}&=\sum_{\ell=1}^d\theta_\ell Y_{k-\ell}Y_{k-d+j}+Y_{k-d+j}Z_k\\
&=\sum_{\substack{\ell=1\\ \ell\neq j-1}}^d\theta_\ell Y_{k-\ell}Y_{k-d+j}+\theta_{j-1}Y_{k-j+1}Y_{k-d+j}+Y_{k-d+j}Z_k
\end{align*}
so that
\begin{displaymath}
\theta_{j-1}Y_{k-j+1}Y_{k-d+j}=Y_kY_{k-d+j}-\sum_{\substack{\ell=1\\ \ell\neq j-1}}^d\theta_\ell Y_{k-\ell}Y_{k-d+j}-
Y_{k-d+j}Z_k
\end{displaymath}
which yields
\begin{align*}
&\theta_d\left(B^{-1}U_kY_k-U_kU_k^T\theta^\ast\right)_j=\theta_dY_kY_{k-j}+\sum_{\substack{\ell=1\\ \ell\neq j-1}}^{d-1}
\theta_\ell Y_{k-j+1}Y_{k-d+1+\ell}-Y_{k-j+1}Y_{k-d+1}\\
&\qquad+Y_kY_{k-d+j}-\sum_{\substack{\ell=1\\ \ell\neq j-1}}^d\theta_\ell Y_{k-\ell}Y_{k-d+j}-
Y_{k-d+j}Z_k\eqdef A_k\left(j\right)-Y_{k-d+j}Z_k
\end{align*}
where
\begin{align*}
A_k\left(j\right)&=\theta_d\left(Y_kY_{k-j}-Y_{k-d}Y_{k-d+j}\right)+\left(Y_kY_{k-d+j}-Y_{k-j+1}Y_{k-d+1}\right)\\
&\qquad+\sum_{\substack{\ell=1\\ \ell\neq j-1}}^{d-1}\theta_\ell\left(Y_{k-j+1}Y_{k-d+1+\ell}-Y_{k-\ell}Y_{k-d+j}\right)\,.
\end{align*}
Hence for $n\geq1$ and $j=2,\ldots,d-1$
\begin{align*}
\sum_{k=1}^nA_k\left(j\right)&=\theta_d\left(\sum_{k=1}^nY_kY_{k-j}-\sum_{k=1}^nY_{k-d}Y_{k-d+j}\right)\\
&\qquad+\left(\sum_{k=1}^nY_kY_{k-d+j}-\sum_{k=1}^nY_{k-j+1}Y_{k-d+1}\right)\\
&\qquad+\sum_{\substack{\ell=1\\ \ell\neq j-1}}^{d-1}\theta_\ell\left(\sum_{k=1}^nY_{k-j+1}Y_{k-d+1+\ell}
-\sum_{k=1}^nY_{k-\ell}Y_{k-d+j}\right)\,.
\end{align*}
Again Lemma~\ref{Lemma:sums} yields $\frac{1}{\sqrt{n}}\sum_{k=1}^nA_k\left(j\right)\rightarrow0$ in probability so that
we get, for $j=2,\ldots,d-1$,
\begin{equation}
\label{Eq:representationcomponentj}
\frac{1}{\sqrt{n}}\sum_{k=1}^n\left(B^{-1}U_kY_k-U_kU_k^T\theta^\ast\right)_j=
-\frac{1}{\theta_d}\frac{1}{\sqrt{n}}\sum_{k=1}^nY_{k-d+j}Z_k+R_n\left(j\right)
\end{equation}
with $R_n\left(j\right)\rightarrow0$ in probability as $n\to\infty$.
\smallskip

\textit{Case 3.} $j=d$. Then by~\eqref{Eq:d}
\begin{displaymath}
\theta_d\left(B^{-1}U_kY_k-U_kU_k^T\theta^\ast\right)_d=
Y_k^2-Y_{k-d+1}^2+\sum_{\ell=1}^{d-1}\theta_\ell
\left(Y_{k-d+1+\ell}Y_{k-d+1}-Y_kY_{k-\ell}\right)\,,
\end{displaymath}
so that
\begin{align*}
&\theta_d\sum_{k=1}^n\left(B^{-1}U_kY_k-U_kU_k^T\theta^\ast\right)_d\\
&\qquad=\sum_{k=1}^n\left(Y_k^2-Y_{k-d+1}^2\right)+\sum_{\ell=1}^{d-1}\theta_\ell\sum_{k=1}^n
\left(Y_{k-d+1+\ell}Y_{k-d+1}-Y_kY_{k-\ell}\right)\,,
\end{align*}
and Lemma~\ref{Lemma:sums} yields
\begin{equation}
\label{Eq:limitcomponentd}
\frac{1}{\sqrt{n}}\sum_{k=1}^n\left(B^{-1}U_kY_k-U_kU_k^T\theta^\ast\right)_d\rightarrow0\quad
\text{in probability as }n\to\infty\,.
\end{equation}

Combining~\eqref{Eq:representationfirstcomponent}, \eqref{Eq:representationcomponentj} and
\eqref{Eq:limitcomponentd} we get
\begin{equation}
\label{Eq:representationIIn}
\frac{1}{\sqrt{n}}\II_n=
\frac{1}{\sqrt{n}}\sum_{k=1}^n\left(B^{-1}U_kY_k-U_kU_k^T\theta^\ast\right)=
\frac{1}{\sqrt{n}}\sum_{k=1}^nh_{\II}\hspace{-2.5pt}\left(U_k\right)Z_k+R_n^{\II}
\end{equation}
with $h_{\II}\hspace{-2.5pt}\left(x\right)=-\frac{1}{\theta_d}\left(x_d,x_{d-1},\ldots,x_2,0\right)^T$, $x\in\RR^d$,
and $R_n^{\II_n}\rightarrow0$ in
probability as $n\to\infty$. From~\eqref{Eq:limitIIIn}, \eqref{Eq:rootIn} and \eqref{Eq:representationIIn} we obtain
\begin{displaymath}
C_n=\frac{1}{\sqrt{n}}\sum_{k=1}^nh\left(U_k\right)Z_k+R_n
\end{displaymath}
with $h\left(x\right)=-\frac{1}{\theta_d}\left(x_d,x_{d-1},\ldots,x_1\right)^T$, $x\in\RR^d$, and
$R_n\rightarrow0$ in probability as $n\to\infty$. Proposition~\ref{Proposition:Clt} implies, because
$E\left(h\left(U_1\right)h\left(U_1\right)^T\right)=\frac{1}{\theta_d^2}\Gamma$,
\begin{equation}
\label{Eq:CltCn}
C_n\stackrel{\D}{\rightarrow}N\left(0,\sigma^2E\left(h\left(U_1\right)h\left(U_1\right)^T\right)\right)
=N\left(0,\frac{\sigma^2}{\theta_d^2}\Gamma\right)\quad\text{as }n\to\infty\,.
\end{equation}
\smallskip

To complete the proof of the proposition we will rely on the following fact: Let $\widetilde{X},\widetilde{X}_n,n\geq1,$ be
$d$-dimensional random vectors defined on some probability space $\left(\Omega,\F,P\right)$ and let $\Omega_n\in\F,n\geq1,$ be events
with $P\left(\Omega_n\right)\rightarrow1$ as $n\to\infty$. Then, as $n\to\infty$,
\begin{equation}
\label{Eq:equivalence}
\widetilde{X}_n\stackrel{\D}{\rightarrow}\widetilde{X}\qquad\text{if and only if}\qquad
\widetilde{X}_n1_{\Omega_n}\stackrel{\D}{\rightarrow}\widetilde{X}\,.
\end{equation}
This follows from the Cram\'{e}r-Slutzky theorem because for all $n\geq1$ and $\varepsilon>0$ we have
$P\left(\left\lVert\widetilde{X}_n-\widetilde{X}_n1_{\Omega_n}\right\rVert\geq\varepsilon\right)\leq P\left(\Omega_n^c\right)\rightarrow0$
as $n\to\infty$ so that $\left\lVert\widetilde{X}_n-\widetilde{X}_n1_{\Omega_n}\right\rVert\rightarrow0$ in probability.
We will apply~\eqref{Eq:equivalence} with
$\Omega_n=\left\{\det\left(\sum_{k-1}^nU_{k-1}U_{k-1}^T\right)>0\right\}$ for $n\geq1$. The sequence 
$\left(\Omega_n\right)_{n\geq1}$ is nondecreasing with $P\left(\Omega_n\right)\rightarrow1$ since
\begin{displaymath}
\det\left(\frac{1}{n}\sum_{k=1}^nU_{k-1}U_{k-1}^T\right)\rightarrow
\det\left(\Gamma\right)\quad\text{almost surely as }n\to\infty
\end{displaymath}
by Lemma~\ref{Lemma:Stronglaws}~(a) because the determinant of a matrix is a continuous function of the matrix components. 
Define $f:\RR^{d\times d}\times\RR^d\rightarrow\RR^d$ by
\begin{displaymath}
f\left(D,x\right)=\left\{
\begin{array}{lcl}
D^{-1}x & , & D\in{\rm GL}\left(d,\RR\right)\,,\\
0       & , & \text{otherwise\,.}
\end{array}\right.
\end{displaymath}
Then $f$ is Borel measurable and continuous at every point in the open subset ${\rm GL}\left(d,\RR\right)\times\RR^d$.
We have on $\Omega_n$
\begin{displaymath}
\sqrt{n}\left(\widehat{\theta}_n-\theta^\ast\right)=f\left(\frac{1}{n}\sum_{k=1}^nU_{k-1}U_{k-1}^T,C_n\right)=:D_n\,.
\end{displaymath}
It follows from Lemma~\ref{Lemma:Stronglaws}~(a) and~\eqref{Eq:CltCn}
\begin{displaymath}
\left(\frac{1}{n}\sum_{k=1}^nU_{k-1}U_{k-1}^T,C_n\right)\stackrel{\D}{\rightarrow}
\left(\Gamma,\frac{\sigma}{\theta_d}\Gamma^{1/2}N_d\right)\,.
\end{displaymath}
The continuous mapping theorem yields
\begin{displaymath}
D_n\stackrel{\D}{\rightarrow}f\left(\Gamma,\frac{\sigma}{\theta_d}\Gamma^{1/2}N_d\right)
=\frac{\sigma}{\theta_d}\Gamma^{-1/2}N_d\,.
\end{displaymath}
This concludes the proof because of~\eqref{Eq:equivalence} .\hfill$\Box$
\medskip

Finally, we comment on the transformation $\theta^\ast=\varphi\left(\theta\right)$, where
\begin{equation}
\label{Eq:definitionphi}
\varphi:\Theta\rightarrow\Theta,\quad\varphi\left(\theta\right)\defeq\left(-\frac{\theta_{d-1}}{\theta_d},\ldots,
-\frac{\theta_1}{\theta_d},\frac{1}{\theta_d}\right)^T
\end{equation}
with $\Theta=\left\{\theta\in\RR^d:\theta_d\neq0\right\}$. Recall that
$\Theta_{pe}=\left\{\theta\in\RR^d:\underline{\rho}\left(B\left(\theta\right)\right)>1\right\}$ and
$\Theta_{s}=\left\{\theta\in\RR^d:\rho\left(B\left(\theta\right)\right)<1\right\}$. 
Then $\Theta_{pe}\subset\Theta$.
\medskip

\begin{Lemma}
\label{Lemma:phiTheta}
$\varphi$ is bijective with $\varphi^{-1}=\varphi$, and we have $\varphi\left(\Theta_{pe}\right)=\Theta_s\cap\Theta$.
\end{Lemma}
\medskip

\textit{Proof.} One easily checks that $\varphi^2\left(\theta\right)=\theta$ for every $\theta\in\Theta$.
Hence, $\varphi$ is  bijective with $\varphi^{-1}=\varphi$. For $\theta\in\Theta$, the eigenvalues of
$B\left(\theta\right)$ are all solutions in $\CC\setminus\left\{0\right\}$ of 
\begin{equation}
\label{Eq:eigenvalueequationBtheta}
\lambda^d-\theta_1\lambda^{d-1}-\cdots-\theta_{d-1}\lambda-\theta_d=0
\end{equation}
while the eigenvalues of $B\left(\varphi\left(\theta\right)\right)$ are all solutions in $\CC\setminus\left\{0\right\}$ of
\begin{equation}
\label{Eq:eigenvalueequationBphitheta}
\lambda^d+\frac{\theta_{d-1}}{\theta_d}\lambda^{d-1}+\cdots+\frac{\theta_1}{\theta_d}\lambda-\frac{1}{\theta_d}=0\,.
\end{equation}
Equation~\eqref{Eq:eigenvalueequationBphitheta} is equivalent to
\begin{displaymath}
\theta_d\lambda^d+\theta_{d-1}\lambda^{d-1}+\cdots+\theta_1\lambda-1=0
\end{displaymath}
which in turn is equivalent to
\begin{displaymath}
\theta_d+\theta_{d-1}\lambda^{-1}+\cdots+\theta_1\lambda^{-\left(d-1\right)}-\lambda^{-d}=0
\end{displaymath}
or, what is the same,
\begin{displaymath}
\lambda^{-d}-\theta_1\lambda^{-\left(d-1\right)}-\cdots-\theta_{d-1}\lambda^{-1}-\theta_d=0\,.
\end{displaymath}
Consequently, $\lambda\in\CC\setminus\left\{0\right\}$ is a solution of~\eqref{Eq:eigenvalueequationBtheta}
if and only if $1/\lambda$ is a solution of~\eqref{Eq:eigenvalueequationBphitheta} so that
\begin{displaymath}
\rho\left(B\left(\varphi\left(\theta\right)\right)\right)=
\frac{1}{\underline{\rho}\left(B\left(\theta\right)\right)}\,.
\end{displaymath}
This implies $\varphi\left(\Theta_{pe}\right)=\Theta_s\cap\Theta$.\hfill$\Box$ 
\medskip

We may extend $\varphi$ to a Borel measurable map on $\RR^d$ by $\varphi\left(x\right)\defeq0\in\RR^d$
if $x\in\RR^d\setminus\Theta$. The ``$\delta$-method'' yields the following central limit theorem
for the consistent estimator $\varphi\left(\widehat{\theta}_n\right)$ of $\theta=\varphi\left(\theta^\ast\right)$.
\medskip

\begin{Corollary}
\label{Corollary:CLTforvarphithetanhat}
We still assume $\theta\in\Theta_{pe}$. We have
\begin{displaymath}
\sqrt{n}\left(\varphi\left(\widehat{\theta}_n\right)-\theta\right)\stackrel{\D}{\rightarrow}
N\left(0,\sigma^2D\Gamma^{-1}D^T\right)\quad\text{as }n\to\infty\,,
\end{displaymath}
where $\Gamma=\Gamma\left(\theta\right)$ and
\begin{displaymath}
D=D\left(\theta\right)=
\begin{pmatrix}
0      & 0 & \ldots & 0 & 0 & 1 & \theta_1     \\
0      & 0 & \ldots & 0 & 1 & 0 & \theta_2     \\
\vdots &   &        &   &   &   & \vdots       \\
1      & 0 & \ldots & 0 & 0 & 0 & \theta_{d-1} \\
0      & 0 & \ldots & 0 & 0 & 0 & \theta_d
\end{pmatrix}
\in\RR^{d\times d}\,.
\end{displaymath}
\end{Corollary}
\medskip

\textit{Proof.} Since $\varphi$ is differentiable at every $x\in\Theta$, it follows from
Proposition~\ref{Proposition:thetanhat} that
\begin{displaymath}
\sqrt{n}\left(\varphi\left(\widehat{\theta}_n\right)-\theta\right)\stackrel{\D}{\rightarrow}
N\left(0,\frac{\sigma^2}{\theta_d^2}\nabla\varphi\left(\theta^\ast\right)\Gamma^{-1}\nabla\varphi\left(\theta^\ast\right)^T\right)
\quad\text{as }n\to\infty\,,
\end{displaymath}
where for $x\in\Theta$
\begin{displaymath}
\nabla\varphi\left(x\right)=
\left(\frac{\partial\varphi_i}{\partial x_j}\left(x\right)\right)_{1\leq i,j\leq d}=
-\frac{1}{x_d}
\begin{pmatrix}
0      & 0 & \ldots & 0 & 0 & 1 & -\frac{x_{d-1}}{x_d} \\
0      & 0 & \ldots & 0 & 1 & 0 & -\frac{x_{d-2}}{x_d} \\
\vdots &   &        &   &   &   & \vdots               \\
1      & 0 & \ldots & 0 & 0 & 0 & -\frac{x_1}{x_d}     \\
0      & 0 & \ldots & 0 & 0 & 0 & \frac{1}{x_d}
\end{pmatrix}
\end{displaymath}
so that $\nabla\varphi\left(\theta^\ast\right)=-\theta_dD$. This yields the assertion.\hfill$\Box$
\medskip

\begin{Remark}\label{Remark:ForwardAR}
A ``forward looking'' variant of the classical autoregressive equation~\eqref{Eq:Yautoregression} is
\begin{equation}\label{Eq:ForwardAR}
Y_n=\theta_1Y_{n+1}+\cdots+\theta_dY_{n+d}+Z_n,\quad n\geq1\,,
\end{equation}
where, as before, $d\in\NN$, $\theta=\left(\theta_1,\ldots,\theta_d\right)^T\in\RR^d$ and
$\left(Z_n\right)_{n\geq1}$ is an i.i.d. sequence of real random variables with $Z_1\in\LL^2\left(P\right)$,
$E\left(Z_1\right)=0$ and $\Var\left(Z_1\right)>0$. For the matrix $B\left(\theta\right)$
and sequence $\left(W_n\right)_{n\geq1}$ as in Section~\ref{Section:Introduction} and $d$-dimensional
vector $V_n$ defined by 
$V_n=\left(Y_n,Y_{n+1},\ldots,Y_{n+d-1}\right)^T$, $n\geq1$, equation~\eqref{Eq:ForwardAR} is equivalent to
\begin{equation}\label{Eq:ForwardARVn}
V_n=B\left(\theta\right)V_{n+1}+W_n,\quad n\geq1\,.
\end{equation}
If $\theta\in\Theta_s$, then the process $\left(V_n\right)_{n\geq1}$
with $V_n=\sum_{k=0}^\infty B\left(\theta\right)^kW_{n+k}$ is stationary (see the proof of 
Lemma~\ref{Lemma:Stationarity}~(i)) and a staightforward computation
shows that it 
satisfies~\eqref{Eq:ForwardARVn} so that the process $\left(Y_n\right)_{n\geq1}$ with 
$Y_n=\pi_1\left(V_n\right)$ for all $n\geq1$ is stationary and satisfies~\eqref{Eq:ForwardAR}.
This means that for given $\theta$ and $\left(Z_n\right)_{n\geq1}$ equation~\eqref{Eq:ForwardAR} has a
stationary solution. This solution is almost surely unique. To see this, let
$\left(\widetilde{Y}_n\right)_{n\geq1}$ be any stationary solution of \eqref{Eq:ForwardAR} and
$\widetilde{V}_n=\left(\widetilde{Y}_n,\widetilde{Y}_{n+1},\ldots,\widetilde{Y}_{n+d-1}\right)$ for $n\geq1$.
By induction,
\begin{displaymath}
\widetilde{V}_n=B\left(\theta\right)^k\widetilde{V}_{n+k}+\sum_{j=0}^{k-1}B\left(\theta\right)^jW_{n+j}\quad
\text{for }n\geq1\text{ and }k\geq0.
\end{displaymath}
Since $\left(\widetilde{Y}_n\right)_{n\geq1}$ is stationary, $\left(\widetilde{V}_n\right)_{n\geq1}$ is bounded in
probability so that, because of $\rho\left(B\left(\theta\right)\right)<1$, for all $n\geq1$,
\begin{displaymath}
B\left(\theta\right)^k\widetilde{V}_{n+k}\rightarrow0\quad\text{in probability as }k\to\infty.
\end{displaymath}
This implies $\widetilde{V}_n=\sum_{j=0}^\infty B\left(\theta\right)^jW_{n+j}=V_n$ almost surely for all
$n\geq1$.
\smallskip

For $\theta_d\neq0$ the solution $\left(Y_n\right)_{n\geq1}$ of \eqref{Eq:ForwardAR}
can also be obtained from the results of Sections~~\ref{Section:Introduction} and~\ref{Section:Features}
by the transformation $\varphi$ and a suitable transformation of the sequence $\left(Z_n\right)_{n\geq1}$
as follows:
\smallskip

Let $\theta\in\Theta_s$ with $\theta_d\neq0$ and the sequence $\left(Z_n\right)_{n\geq1}$ 
with $Z_1\in\LL^2\left(P\right)$, $E\left(Z_1\right)=0$ and $\Var\left(Z_1\right)>0$ be given. Then
$\theta\in\Theta_s\cap\Theta$ and $\theta^*=\left(\theta^*_1,\ldots,\theta^*_d\right)^T
\defeq\varphi\left(\theta\right)\in\Theta_{pe}$ by Lemma~\ref{Lemma:phiTheta}. The sequence
$\left(Z^*_{n}\right)_{n\geq1}$ with $Z^*_n\defeq-\frac{1}{\theta_d}Z_n$ for $n\geq1$ is
i.i.d. with $\Var\left(Z^*_1\right)>0$ and $E\left(Z^*_1\right)=0$. Therefore, by the
results of Sections~\ref{Section:Introduction} and~\ref{Section:Features}, there exists a stationary process
$\left(Y^*_{n-d+1}\right)_{n\geq0}$ with
\begin{equation*}
Y^*_n=\theta^*_1Y^*_{n-1}+\cdots+\theta^*_dY^*_{n-d}+
Z^*_n\quad\text{for all }n\geq1
\end{equation*}
which means
\begin{equation*}
Y^*_n=-\frac{\theta_{d-1}}{\theta_d}Y^*_{n-1}-\cdots-
\frac{\theta_1}{\theta_d}Y^*_{n-d+1}+\frac{1}{\theta_d}Y^*_{n-d}-
\frac{1}{\theta_d}{Z}_n\quad\text{for all }n\geq1
\end{equation*}
which in turn is equivalent to
\begin{equation*}
Y^*_{n-d}=\theta_1Y^*_{n-d+1}+\cdots+\theta_{d-1}Y^*_{n-1}+\theta_dY^*_n+Z_n
\quad\text{for all }n\geq1\,.
\end{equation*}
Therefore, the process $\left(Y_n\right)_{n\geq1}$ with $Y_n\defeq Y^*_{n-d}$ for $n\geq1$ is
stationary with
\begin{equation*}
Y_n=\theta_1Y_{n+1}+\cdots+\theta_{d-1}Y_{n+d-1}+\theta_dY_{n+d}+Z_n\quad\text{for all }n\geq1\,,
\end{equation*}
that is, $\left(Y_n\right)_{n\geq1}$ is a stationary solution of~\eqref{Eq:ForwardAR} for the given
$\theta$ and $\left(Z_n\right)_{n\geq1}$. Clearly, one can
go the other way round and construct a stationary solution of~\eqref{Eq:Yautoregression} from a
stationary solution of~\eqref{Eq:ForwardAR} provided that $\theta_d\neq0$. Therefore, the probabilistic
models defined by \eqref{Eq:Yautoregression} for $\theta\in\Theta_{pe}$ and by \eqref{Eq:ForwardAR} for
$\theta\in\Theta_s\cap\Theta$ are essentially the same.
\smallskip

Note that $Y^*_{n-d+1}=\pi_d\left(U^*_n\right)$ for $U^*_n=-\sum_{k=1}^\infty B\left(\theta^*\right)^{-k}W^*_{n+k}$
with $W^*_n=\left(Z^*_n,0,\right.$ $\left.\ldots,0\right)^T$ for $n\geq1$; see \eqref{Eq:RepresentationUn}. 
On the other hand, by almost sure
uniqueness of the stationary solution of \eqref{Eq:ForwardAR}, $Y^*_{n-d}=\pi_1\left(V_n\right)$ almost surely so
that for $\theta\in\Theta_s\cap\Theta$ and $n\geq1$ almost surely
\begin{align*}
&\pi_1\left(\sum_{k=0}^\infty B\left(\theta\right)^kW_{n+k}\right)=\pi_1\left(V_n\right)=Y^*_{n-d}
=\pi_d\left(U^*_{n-1}\right)\\
&\qquad=\pi_d\left(-\sum_{k=1}^\infty B\left(\theta^*\right)^{-k}W^*_{n-1+k}\right)
=\pi_d\left(\sum_{k=1}^\infty B\left(\theta^*\right)^{-k}\frac{1}{\theta_d}W_{n-1+k}\right)\\
&\qquad=\pi_d\left(\sum_{k=0}^\infty B\left(\varphi\left(\theta\right)\right)^{-k-1}\frac{1}{\theta_d}W_{n+k}\right)\,.
\end{align*}
For $d=1$ the resulting equation is obvious because $B\left(\theta\right)=\theta_1$, $\varphi\left(\theta\right)=\frac{1}{\theta_1}$
and hence $B\left(\varphi\left(\theta\right)\right)^{-1}=B\left(\frac{1}{\theta_1}\right)^{-1}=\left(\frac{1}{\theta_1}\right)^{-1}=\theta_1$
so that $B\left(\theta\right)^k=B\left(\varphi\left(\theta\right)\right)^{-k-1}\frac{1}{\theta_1}$ for all $k\geq0$.
For $d\geq2$, however, the equation does not seem to be so obvious.
\smallskip

For $d=1$ and distributions of $Z_1$ with heavy tails the forward looking autoregressive model~\eqref{Eq:ForwardAR}
has been studied in~\cite{GourierouxZakoian}.
For arbitrary $d\in\NN$ the model defined by~\eqref{Eq:ForwardAR} is a special case of the model considered
in~\cite{LanneSaikkonen}.
\hfill$\Box$
\end{Remark}

\renewcommand\refname

\end{document}